# Gromov's measure equivalence and rigidity of higher rank lattices

By Alex Furman

## Abstract

In this paper the notion of Measure Equivalence (ME) of countable groups is studied. ME was introduced by Gromov as a measure-theoretic analog of quasi-isometries. All lattices in the same locally compact group are Measure Equivalent; this is one of the motivations for this notion. The main result of this paper is ME rigidity of higher rank lattices: any countable group which is ME to a lattice in a simple Lie group $G$ of higher rank, is commensurable to a lattice in $G$.

## 1. Introduction and statement of main results

This is the first in a sequence of two papers on rigidity aspects of measure-preserving group actions, the second being [Fu]. In this paper we discuss the following equivalence relation between (countable) groups, which was introduced by Gromov:

*Definition* 1.1 ([Gr, 0.5.E]). Two countable groups $\Gamma$ and $\Lambda$ are said to be *Measure Equivalent* (ME) if there exist commuting, measure-preserving, free actions of $\Gamma$ and $\Lambda$ on some infinite Lebesgue measure space $(\Omega, m)$, such that the action of each of the groups $\Gamma$ and $\Lambda$ admits a finite measure fundamental domain. The space $(\Omega, m)$ with the actions of $\Gamma$ and $\Lambda$ will be called a ME *coupling* of $\Gamma$ with $\Lambda$.

The basic example of ME groups are lattices in the same locally compact group:

*Example* 1.2. Let $\Gamma$ and $\Lambda$ be lattices in the same locally compact second countable (lcsc) group $G$. Since $G$ contains lattices it is necessarily unimodular, so its Haar measure $m_G$ is invariant under the left $\Gamma$-action $\gamma : g \mapsto \gamma^{-1} g$, and the right $\Lambda$-action $\lambda : g \mapsto g\,\lambda$, which obviously commute. Hence $(G, m_G)$ forms a ME coupling for $\Gamma$ and $\Lambda$.

The relation of ME between groups can be considered as a measure-theoretic analog of quasi-isometry between groups, due to the following beautiful observation of Gromov:

THEOREM ([Gr, 0.2.C]). *Two finitely generated groups $\Gamma$ and $\Lambda$ are quasi-isometric if and only if they are* topologically equivalent *in the following sense: there exist commuting, continuous actions of $\Gamma$ and $\Lambda$ on some locally compact space $X$, such that the actions of each of the groups is properly discontinuous and cocompact.*

A typical example of such topological equivalence consists of a locally compact group and any two of its uniform (i.e. cocompact) lattices which act by translations from the left and from the right. Thus all uniform lattices in the same lcsc group are topologically equivalent (i.e. quasi-isometric); however typically nonuniform lattices are not quasi-isometric to the uniform ones.

Examples of lattices suggest that ME is a (strictly) weaker relation than topological equivalence (i.e. quasi-isometry). For general groups this is not known; quasi-isometric groups admit topological coupling (in the above sense), but it is not clear whether one can always find such a coupling carrying an invariant measure.

The notion of ME is also closely related to (Weak) Orbit Equivalence of measure-preserving free actions of groups on probability spaces (c.f. [Fu]). This connection was observed by R. Zimmer and was apparently known also to M. Gromov. More precisely, one has

THEOREM (see [Fu, 3.2 and 3.3]). *Let $\Gamma$ and $\Lambda$ be two countable groups which admit Weakly Orbit Equivalent, essentially free, measure-preserving actions on Lebesgue probability spaces $(X, \mu)$ and $(Y, \nu)$. Then $\Gamma$ and $\Lambda$ are measure equivalent. Moreover, there exists a ME coupling $(\Omega, m)$ of $\Gamma$ with $\Lambda$, such that the $\Gamma$-action on the quotient $\Omega/\Lambda$ and the $\Lambda$-action on $\Gamma \backslash \Omega$ are isomorphic to $(X, \Gamma)$ and $(Y, \Lambda)$ respectively.*

*Conversely, given an ME coupling $(\Omega, m)$ of $\Gamma$ with $\Lambda$, the left $\Gamma$-action on $\Omega/\Lambda$ and the right $\Lambda$-action on $\Gamma \backslash \Omega$ give weakly orbit equivalent measure-preserving actions on probability spaces.*

The concept of (weak) Orbit Equivalence has been thoroughly studied in the framework of ergodic theory. In particular, Ornstein and Weiss [OW] have proved (generalizing previous work of Dye; see also [CFW]) that all (not necessarily free) ergodic measure-preserving actions of all countable amenable groups are Orbit Equivalent. Hence all countable amenable groups are ME. On the other hand, it is well known (c.f. [Zi3, 4.3.3]) that nonamenable groups are not ME to amenable ones. Hence

COROLLARY 1.3. *The* ME *class of integers $\mathbb{Z}$ consists precisely of all countable amenable groups.*

This fact shows that ME is a rather weak relation, compared to quasi-isometries. For example, the class of groups which are topologically equivalent (i.e. quasi-isometric) to $\mathbb{Z}$ consists only of finite extensions of $\mathbb{Z}$. By considering various amenable groups, one also observes that many quasi-isometric invariants, such as growth functions, finite generation, cohomological dimension etc. are not preserved under a ME relation. The property of being a word hyperbolic group is not preserved by ME either, because $\mathrm{SL}_2(\mathbb{C})$ contains both word hyperbolic and not-word hyperbolic lattices.

However, ME is convenient for transferring unitary representation invariants. Given a ME coupling $(\Omega, m)$ of $\Gamma$ with $\Lambda$ one can induce unitary $\Lambda$-representations to unitary $\Gamma$-representations (see Section 8), which gives

COROLLARY 1.4. *Kazhdan's property T is a* ME *invariant.*

We remark, that unlike amenability, it is still an open problem whether property T is a quasi-isometric invariant.

Lattices in (semi-) simple Lie groups form an especially interesting class of examples of discrete groups. The program of quasi-isometric classification of lattices has been recently completed by a sequence of works by Eskin, Farb, Kleiner, Leeb and Schwartz (see [Sc], [FS], [EF], [Es], [KL], and the surveys [GP] and [Fa]). It is known now, that any group which is quasi-isometric to a lattice in a semisimple Lie group is commensurable to a lattice in the same Lie group, while lattices in the same simple Lie group split into several quasi-isometric classes with the class of uniform lattices being one of them.

Consider a similar classification program in the ME context. Recall that all lattices in the same group are automatically ME (Example 1.2). The correspondence between ME and (weak) Orbit Equivalence of measure-preserving actions, mentioned above, enables us to translate Zimmer's work on Orbit Equivalence (which is based on superrigidity for measurable cocycles theorem [Zi1]), to the ME framework. This almost direct translation shows that a lattice $\Gamma$ in a higher rank simple Lie group is not ME to a lattice in any other semisimple Lie group. A further result [Zi2] of Zimmer in this direction shows that if a lattice $\Gamma$ in a higher rank simple Lie group $G$ is ME to a countable group $\Lambda$, which is just known to admit a *linear representation* with an infinite image, then $\Lambda$ is commensurable to a lattice in $\mathrm{Ad}\,G$.

In the present paper we prove that any group $\Lambda$, which is ME to $\Gamma$, is essentially linear, removing completely the linearity assumptions in the above results. Hence for a higher rank simple Lie group $G$, the collection of all its lattices (up to finite groups) forms a single ME class. Moreover, we show that any ME coupling of $\Gamma$ with $\Lambda$ has a standard ME coupling as a factor. More precisely

THEOREM (3.1). *Let $\Gamma$ be a lattice in a simple, connected Lie group $G$ with finite center and $\mathbb{R}-\mathrm{rk}(G) \geq 2$. Let $\Lambda$ be an arbitrary countable group, which is Measure Equivalent to $\Gamma$. Then there exists a finite index subgroup $\Lambda' \subseteq \Lambda$ and an exact sequence*

$$1 \longrightarrow \Lambda_0 \longrightarrow \Lambda' \overset{\rho}{\longrightarrow} \Lambda_1 \longrightarrow 1$$

*where $\Lambda_0 \lhd \Lambda$ is finite and $\Lambda_1 \subset \mathrm{Ad}\, G$ is a lattice. Moreover, if $(\Omega, m)$ is a ME coupling of $\Gamma$ with $\Lambda$, then there exists a unique measurable map $\Phi : \Omega \to \mathrm{Ad}\, G$ such that*

$$\Phi(\gamma\, \omega\, \lambda') = \mathrm{Ad}(\gamma)\, \Phi(x)\, \rho(\lambda'), \qquad (\gamma \in \Gamma, \ \lambda' \in \Lambda').$$

*The measure $\Phi_* m$ is a convex combination of an atomic measure and the Haar measure on $\mathrm{Ad}\, G$; disintegration of $m$ with respect to $\Phi_* m$ consists of probability measures.*

The main substance of the theorem is the construction of a virtually faithful representation for the unknown group $\Lambda$, using just the fact that $\Lambda$ is ME to a higher rank lattice $\Gamma$.

*Open questions.*

1. What is the complete ME classification of lattices in (semi)simple Lie groups? By Theorem 3.1 it remains to classify lattices in rank one simple Lie groups. Since property T is a ME invariant, lattices in $\mathrm{Sp}(n, 1)$ and $F_4$ are not ME to lattices in $\mathrm{SO}(n, 1)$ or in $\mathrm{SU}(n, 1)$. Lattices in different $\mathrm{Sp}(n, 1)$ are not ME to each other; this follows from the work of Cowling and Zimmer [CZ] which uses von Neumann algebra invariants. So the main question is whether ME distinguishes lattices in different groups among $\mathrm{SO}(n, 1)$ and $\mathrm{SU}(n, 1)$.

2. Characterize or describe general properties of the class of (countable) groups which are ME to a free group. Note that this class contains (all) lattices in $\mathrm{SL}_2(\mathbb{R})$, $\mathrm{SL}_2(\mathbb{Q}_p)$ and in $\mathrm{Aut}(T)$ – the group of automorphisms of a regular tree.

3. Find other ME invariants, besides amenability and property T.

4. Is it true that any two quasi-isometric groups are ME? Note that if true, this, combined with Corollary 1.4, would imply that Kazhdan's property T is a quasi-isometric invariant.

*Acknowledgments.* I am deeply grateful to Robert Zimmer, Benson Farb and Alex Eskin for their support, constant encouragement and for many enlightening conversation on quasi-isometries, orbit equivalence and related topics.

## 2. Definitions, notations and some basic properties

Inspired by the fact that any two lattices in the same group are ME (Example 1.2), we shall use a similar left and right notation for a general ME coupling. More precisely, given a ME coupling $(\Sigma, \sigma)$ of arbitrary countable groups $\Gamma$ and $\Lambda$ we shall denote the actions by:

$$\gamma : x \mapsto \gamma^{-1} x \qquad \lambda : x \mapsto x\,\lambda \qquad (x \in \Sigma,\ \gamma \in \Gamma,\ \lambda \in \Lambda).$$

Thus, by saying that $(\Sigma, \sigma)$ is a ME coupling of $\Gamma$ with $\Lambda$ we shall mean that $\Gamma$ acts on $(\Sigma, \sigma)$ from the left and $\Lambda$ acts from the right. Thus, formally, the definition of ME becomes asymmetric.

*Duals and compositions of* ME *couplings.* Using the terminology of Measure Equivalence we should check that the relation is indeed an equivalence relation. Showing this we shall establish some notation and terminology to be used in the sequel.

*Reflexivity.* Any countable group $\Gamma$ is ME to itself: consider the left and the right $\Gamma$ action on itself with the Haar (counting) measure $m_\Gamma$.

*Symmetry.* Given a ME coupling $(\Sigma, \sigma)$ of $\Gamma$ with $\Lambda$, one can consider the *dual* ME coupling $(\check{\Sigma}, \check{\sigma})$ of $\Lambda$ with $\Gamma$, defined formally as follows: the measure spaces $(\Sigma, \sigma)$ and $(\check{\Sigma}, \check{\sigma})$ are isomorphic with $x \mapsto \check{x}$ being the isomorphism; the left $\Lambda$-action and the right $\Gamma$-actions on $\check{\Sigma}$ are defined by $\lambda^{-1}\,\check{x} = (x\,\lambda)\check{\ }$ and $\check{x}\,\gamma = (\gamma^{-1}\,x)\check{\ }$.

*Transitivity.* Given a ME coupling $(\Sigma, \sigma)$ of $\Gamma$ with $\Lambda$, and a ME coupling $(\Sigma', \sigma')$ of $\Lambda$ with $\Delta$, define the *composition* ME coupling of $\Gamma$ with $\Delta$, denoted by $(\Sigma \times_\Lambda \Sigma',\ \sigma \times_\Lambda \sigma')$. The latter is constructed as follows: consider the commuting measure-preserving actions of $\Gamma$, $\Lambda$ and $\Delta$ on the product space $(\Sigma \times \Sigma',\ \sigma \times \sigma')$:

$$(2.1) \qquad \begin{aligned} \gamma &: (x, y) \mapsto (\gamma^{-1}x, y), \quad \delta : (x, y) \mapsto (x, y\,\delta), \\ \lambda &: (x, y) \mapsto (x\,\lambda,\ \lambda^{-1}y) \end{aligned}$$

for $x \in \Sigma$, $y \in \Sigma'$ and $\gamma \in \Gamma$, $\lambda \in \Lambda$, $\delta \in \Delta$. The composition ME coupling is the action of $\Gamma$ and $\Delta$ on the space of the $\Lambda$-orbits in $\Sigma \times \Sigma'$, equipped with the natural factor measure. This measure becomes apparent when the space $\Sigma \times_\Lambda \Sigma'$ of $\Lambda$-orbits is identified with the spaces

$$\Sigma \times_\Lambda \Sigma' \cong \Sigma \times X' \cong Y \times \Sigma'$$

where $Y \subset \Sigma$ and $X' \subset \Sigma'$ are some fundamental domains for $\Sigma/\Lambda$ and $\Lambda \backslash \Sigma'$, respectively. From these identifications one readily sees that $\Gamma$ and $\Delta$ admit finite measure fundamental domains, namely $X \times X'$ for $\Gamma \backslash (\Sigma \times X')$ and $Y \times Y'$ for $(Y \times \Sigma')/\Delta$, where $X \subset \Sigma$ and $Y' \subset \Sigma'$ are (left) $\Gamma$- and (right) $\Delta$-fundamental domains.

*Example* 2.1.  To clarify the above terminology consider our illustrative Example 1.2: let $G$ be a lcsc group with the Haar measure $m_G$, $\Gamma$ and $\Lambda$ lattices in $G$, and take the ME coupling $(G, m_G)$ of $\Gamma$ with $\Lambda$. The dual ME coupling $(\check{G}, \check{m}_G)$ of $\Lambda$ with $\Gamma$ can be identified also with the natural ME coupling $(G, m_G)$. This time one considers the left $\Lambda$-action and the right $\Gamma$-action, the identification of $(\check{G}, \check{m}_G)$ with $(G, m_G)$ being given by $\check{x} = x^{-1}$.

Let $\Gamma_1$, $\Gamma_2$ and $\Gamma_3$ be three lattices in the same lcsc group $G$. Considering the composition $G \times_{\Gamma_2} G$ of the natural ME couplings of $\Gamma_1$ with $\Gamma_2$ and $\Gamma_2$ with $\Gamma_3$, one can check that this ME coupling of $\Gamma_1$ with $\Gamma_3$ admits a $\Gamma_1 \times \Gamma_3$-equivariant map onto the natural ME coupling by $G$. Indeed consider the map $G \times G \to G$ given by $(x, y) \mapsto xy$ (the latter is multiplication in $G$) which factors through the space of $\Gamma_2$-orbits $G \times_{\Gamma_2} G$. The corresponding map $G \times_{\Gamma_2} G \to G$ is the required map, with $G/\Gamma_2$ (or $\Gamma_2 \backslash G$) being the fiber.

*Cocycles $\alpha : \Gamma \times \Sigma/\Lambda \to \Lambda$ and $\beta : \Gamma \backslash \Sigma \times \Lambda \to \Gamma$.* Let $(\Sigma, \sigma)$ be a ME coupling of two countable groups $\Gamma$ and $\Lambda$ and let $(X, \mu)$ and $(Y, \nu)$ be some fundamental domains for $\Sigma/\Lambda$ and $\Gamma \backslash \Sigma$, respectively (here $\mu = \sigma|_X$ and $\nu = \sigma|_Y$). Define a measurable function $\alpha : \Gamma \times \Sigma \to \Lambda$ by the following rule: given $x \in X$ and $\gamma \in \Gamma$, let $\alpha(\gamma, x)$ be the unique element $\lambda \in \Lambda$, satisfying $\gamma\, x \in X\, \lambda$. In a similar way, given a fundamental domain $Y$ for $\Gamma \backslash \Sigma$, we define $\beta : Y \times \Lambda \to \Gamma$ by the following: given $y \in Y$ and $\lambda \in \Lambda$, set $\beta(y, \lambda) = \gamma \in \Gamma$ if $y\, \lambda \in \gamma\, Y$.

With these definitions the natural actions of $\Gamma$ on $\Sigma/\Lambda \cong X$ and $\Lambda$ on $\Gamma \backslash \Sigma \cong Y$ can be described by the formulas

$$\begin{aligned}
\gamma \cdot x &= \gamma\, x\, \alpha(\gamma, x)^{-1} \in X & (x \in X,\ \gamma \in \Gamma) \\
y \cdot \lambda &= \beta(y, \lambda)^{-1}\, y\, \lambda \in Y & (y \in Y,\ \lambda \in \Lambda)
\end{aligned}$$

where the actions in the right-hand sides are the original $\Gamma$ and $\Lambda$ actions on $\Sigma$, while the actions on the left-hand sides (denoted by the dot) represent the natural actions on the spaces of orbits $\Sigma/\Lambda$ and $\Gamma \backslash \Sigma$.

From the definitions it follows directly that $\alpha$ and $\beta$ are (left and right) cocycles; i.e.

$$\begin{aligned}
\alpha(\gamma_1\, \gamma_2, x) &= \alpha(\gamma_1, \gamma_2 \cdot x)\, \alpha(\gamma_2, x) \\
\beta(y, \lambda_1\, \lambda_2) &= \beta(y, \lambda_1)\, \beta(y \cdot \lambda_1, \lambda_2).
\end{aligned}$$

Moreover, choosing another fundamental domain, say $X'$ for $\Sigma/\Lambda$, would result in a cocycle $\alpha'$ which is measurably cohomologous to $\alpha$. More precisely, if $\theta : X \to X'$ is the isomorphism, then

$$\alpha'(\gamma, \theta(x)) = \xi(\gamma \cdot x)^{-1}\, \alpha(\gamma, x)\, \xi(x)$$

where $\xi : X \to \Lambda$ is chosen so that $\theta(x) = x\, \xi(x) \in X'$. Similar statements hold for $\beta$.

Having this in mind, we can talk about the "canonical" cocycles $\alpha : \Gamma \times \Sigma/\Lambda \to \Lambda$ and $\beta : \Gamma \backslash \Sigma \times \Lambda \to \Gamma$ meaning the corresponding cohomological class of measurable cocycles.

*Ergodicity of* ME *couplings.* A ME coupling $(\Sigma, \sigma)$ of $\Gamma$ and $\Lambda$ is said to be *ergodic* if the left $\Gamma$-action on $\Sigma/\Lambda$ and the right $\Lambda$-action on $\Gamma \backslash \Sigma$ are ergodic.

LEMMA 2.2. *Let* $(\Sigma, \sigma)$ *be a* ME *coupling of two countable groups* $\Gamma$ *and* $\Lambda$.

1. *The* $\Gamma$-*action on* $\Sigma/\Lambda$ *is ergodic if and only if the* $\Lambda$-*action on* $\Gamma \backslash \Sigma$ *is ergodic.*

2. *The* $\Gamma \times \Lambda$-*invariant measure* $\sigma$ *on* $\Sigma$ *can be disintegrated in the form* $\sigma = \int \sigma_t \, d\eta(t)$, *where* $\eta$ *is some probability measure, such that for* $\eta$-*almost every* $t$ *the measure* $\sigma_t$ *is* $\Gamma \times \Lambda$-*invariant and* $(\Sigma, \sigma_t)$ *forms an* ergodic ME *coupling of* $\Gamma$ *with* $\Lambda$.

3. *The composition coupling* $(\Sigma \times_\Lambda \check{\Sigma}, \sigma \times_\Lambda \check{\sigma})$ *of* $\Gamma$ *with* $\Gamma$ *is ergodic if and only if the* $\Lambda$-*action on* $\Gamma \backslash \Sigma$ *is weakly mixing, in which case any composition* ME *coupling* $(\Sigma \times_\Lambda \Sigma', \sigma \times_\Lambda \sigma')$ *of* $\Gamma$ *with* $\Delta$ *is ergodic, provided that the* ME *coupling* $(\Sigma', \sigma')$ *of* $\Lambda$ *with* $\Delta$ *is ergodic.*

*Proof.* Let $(X, \mu)$ and $(Y, \nu)$ be some fundamental domains in $(\Sigma, \sigma)$ with respect to $\Lambda$ and $\Gamma$ actions.

1. Note that $\Gamma$ is ergodic on $(X, \mu)$ if and only if $\Gamma \times \Lambda$ is ergodic on $(\Sigma, \sigma)$, which happens if and only if $\Lambda$ is ergodic on $(Y, \nu)$.

2. Consider the $\Gamma$ on $(X, \mu)$. Using ergodic decomposition, write $\mu = \int \mu_t \, d\eta(t)$, where $\eta$ is some probability measure and $\eta$-a.e. $\mu_t$ are $\Gamma$-invariant and ergodic *probability* measures on $X$. Let $\sigma_t$ be the lifting of $\mu_t$ from $X \cong \Sigma/\Lambda$ to $\Sigma$. We have $\sigma = \int \sigma_t \, d\eta(t)$ with $\sigma_t$ being $\Gamma \times \Lambda$-invariant. Let $\nu_t = \sigma_t|_Y$, then

$$\nu(Y) = \int \nu_t(Y) \, d\eta(t) < \infty.$$

Hence for $\eta$-a.e. $t$ the measure $\nu_t$ is *finite*. Therefore, for $\eta$-a.e. $t$, $(\Sigma, \sigma_t)$ is a ME coupling of $\Gamma$ with $\Lambda$, which is ergodic by 1.

3. Let $(\Sigma', \sigma')$ be some ME coupling of $\Lambda$ with $\Delta$, and let $(Z, \eta)$ be a fundamental domain for the right $\Delta$-action. Observe that the composition ME coupling $(\Sigma \times_\Lambda \Sigma', \sigma \times_\Lambda \sigma')$ of $\Gamma$ with $\Delta$ is ergodic if and only if the action (2.1) of $\Gamma \times \Lambda \times \Delta$ on $(\Sigma \times \Sigma', \sigma \times \sigma')$ is ergodic. The latter happens if and only if the $\Lambda$-action on $(Y \times Z, \nu \times \eta)$, given by

$$\lambda : (y, z) \mapsto (y \cdot \lambda, \lambda^{-1} \cdot z) \ ,$$

is ergodic. This is just the product (or diagonal) action of $\Lambda$ on $(Y, \nu) \times (Z, \eta)$. The assertion follows from the standard facts on weakly mixing group actions (c.f. [BR]). $\qquad \square$

## 3. The main result and an outline of its proof

We shall state our main result in a slightly different form, using the group of all automorphisms $\operatorname{Aut}(\operatorname{Ad} G)$ of $\operatorname{Ad} G$. The relationship between $G$, its center $Z(G)$, and the groups $\operatorname{Inn}(G) = \operatorname{Ad} G$ and $\operatorname{Aut}(\operatorname{Ad} G)$ is described by the exact sequence

$$1 \longrightarrow Z(G) \longrightarrow G \xrightarrow{\operatorname{Ad}} \operatorname{Ad} G \longrightarrow \operatorname{Aut}(\operatorname{Ad} G) \xrightarrow{o} \operatorname{Out}(\operatorname{Ad} G) \longrightarrow 1$$

where $\operatorname{Ad}(g) : h \mapsto g^{-1} \, h \, g$. Note that in our case both $\operatorname{Out}(\operatorname{Ad} G)$ and $Z(G)$ are finite groups.

THEOREM 3.1 (Measure Equivalence rigidity for higher rank lattices). *Let $G$ be a simple connected Lie group with a finite center and $\mathbb{R}\!-\!\operatorname{rk}(G) \geq 2$. Let $\Gamma \subset G$ be a lattice and let $\Lambda$ be an* arbitrary *countable group, which is Measure Equivalent to $\Gamma$. Then there exists a homomorphism $\rho : \Lambda \to \operatorname{Aut}(\operatorname{Ad} G)$ with a finite kernel $\Lambda_0 = \operatorname{Ker}(\rho)$ and the image $\operatorname{Im}(\rho) = \Lambda_1 \subset \operatorname{Aut}(\operatorname{Ad} G)$ being a lattice in $\operatorname{Aut}(\operatorname{Ad} G)$. Moreover, if $(\Sigma, \sigma)$ is the ME coupling between $\Gamma$ and $\Lambda$, then there exists a unique measurable map $\Phi : \Sigma \to \operatorname{Aut}(\operatorname{Ad} G)$ satisfying*

$$\Phi(\gamma^{-1} \, x \, \lambda) = \operatorname{Ad}(\gamma)^{-1} \, \Phi(x) \, \rho(\lambda)$$

*for $\gamma \in \Gamma$, $\lambda \in \Lambda$ and $\sigma$-a.e. $x \in \Sigma$. The projection $\Phi_*\sigma$ of the measure $\sigma$ to $\operatorname{Aut}(\operatorname{Ad} G)$ is a convex combination of an atomic measure and Haar measures on $\operatorname{Ad} G$-cosets in $\operatorname{Aut}(\operatorname{Ad} G)$, where the fibers of the disintegration of $\sigma$ with respect to $\Phi_*\sigma$ are probability measures.*

In the statement of the theorem in the introduction we take

$$\Lambda' = \operatorname{Ker}(o \circ \rho : \Lambda \to \operatorname{Out}(\operatorname{Ad} G)).$$

*Outline of the proof.* The proof of the main theorem contains four steps, described briefly below.

*Step* 1. *Analysis of self* ME *couplings of* $\Gamma$. We consider a general ME coupling $(\Omega, m)$ of $\Gamma$ with itself and show that it has a uniquely defined $\Gamma \times \Gamma$-equivariant measurable mapping $\Phi : \Omega \to \operatorname{Aut}(\operatorname{Ad} G)$. The main tool in the construction of $\Phi$ is Zimmer's superrigidity for cocycles. Uniqueness of $\Phi$ is proved by an argument on smoothness of algebraic actions. Ratner's theorem is used to identify the image $\Phi_*m$ of the measure $m$.

*Step* 2. *Construction of the representation.* Given a ME coupling $(\Sigma, \sigma)$ of $\Gamma$ with an unknown group $\Lambda$ we construct a representation $\rho : \Lambda \to \operatorname{Aut}(\operatorname{Ad} G)$. The idea of the construction is to consider the composition ME coupling $\Omega = (\Sigma \times_\Lambda \Lambda \times_\Lambda \check{\Sigma})$ of $\Gamma$ with $\Gamma$, and to use the factoring map $\Phi : \Omega \to \operatorname{Aut}(\operatorname{Ad} G)$. The main point of the proof is to show that for a.e. fixed $(x, \check{y}) \in \Sigma \times \check{\Sigma}$ certain expression in terms of $\Phi([x, \lambda, \check{y}])$ gives a representation $\rho_{x,y} : \Lambda \to \mathrm{Aut}(\mathrm{Ad}\, G)$. It turns out, that the representations $\rho_{x,y} = \rho_x$ do not depend on $y$, and that different values of $x$ give rise to equivalent representations.

*Step* 3. *Finiteness of the kernel.* The construction of the representations $\rho_x$ enables to show that the common kernel $\Lambda_0 = \mathrm{Ker}(\rho_x : \Lambda \to \mathrm{Aut}(\mathrm{Ad}\, G))$ is finite.

*Step* 4. *The image is a lattice.* Once obtained, the linear representation $\rho$ enables us to apply Zimmer's result [Zi2], stating that in this case $\rho(\Lambda)$ forms a lattice in $\mathrm{Aut}(\mathrm{Ad}\, G)$. This argument relies on another application of superrigidity for cocycles (with real and $p$-adic targets).

*Remark* 3.2. A reader familiar with the proofs of quasi-isometric rigidity results may recognize some lines of similarity in the scheme of the proof:

- For any finitely generated group $\Gamma$, there exists an associated group $\mathrm{QI}(\Gamma)$ of its self-quasi-isometries, which is extremely useful in the study of quasi-isometric properties of $\Gamma$. In particular, given any quasi-isometry $q : \Lambda \to \Gamma$ of an unknown group $\Lambda$ to a known $\Gamma$, one automatically obtains a representation

$$\rho_q : \Lambda \xrightarrow{\tau} \mathrm{QI}(\Lambda) \xrightarrow{q_*} \mathrm{QI}(\Gamma) \qquad \rho_q(\lambda) = q_* \circ \tau(\lambda) \circ q_*^{-1},$$

where $\tau$ is the representation by translations, and $q_* : \mathrm{QI}(\Lambda) \cong \mathrm{QI}(\Gamma)$ is the isomorphism corresponding to $q$. Hence studying quasi-isometric rigidity amounts to the identification of the group $\mathrm{QI}(\Gamma)$, analysis of the image $\rho_q(\Lambda) \subseteq \mathrm{QI}(\Gamma)$, and a proof of the finiteness of the kernel $\mathrm{Ker}(\rho_q)$.

- In a somewhat analogous framework of ME we do not see a reasonable general definition of a ME analog of the group $\mathrm{QI}(\Gamma)$, and therefore do not have an abstract construction of a representation as above. Nevertheless, Step 1 of our proof can be interpreted as an identification of a (nonexisting) ME analog of the group $\mathrm{QI}(\Gamma)$ with $\mathrm{Aut}(\mathrm{Ad}\, G)$. Step 2 of the proof traces the construction of $\rho_q$ with most of the effort devoted to the proof that the construction indeed gives a representation.

## 4. Self ME couplings of higher rank lattices

THEOREM 4.1. *Let $G$ be as in Theorem 3.1. Let $\Gamma_1$, $\Gamma_2 \subset G$ be lattices and let $(\Omega, m)$ be a ME coupling of $\Gamma_1$ with $\Gamma_2$. Then there exists a unique measurable map $\Phi_\Omega : \Omega \to \mathrm{Aut}(\mathrm{Ad}\, G)$, satisfying*

$$(4.1) \qquad \Phi_\Omega(\gamma_1\, \omega\, \gamma_2) = \mathrm{Ad}(\gamma_1)\, \Phi_\Omega(\omega)\, \mathrm{Ad}(\gamma_2)$$

*for $m$-a.e. $\omega \in \Omega$ and all $\gamma_1 \in \Gamma_1$ and $\gamma_2 \in \Gamma_2$. The projection $\Phi_{\Omega*}m$ of $m$ is a convex combination of an atomic and Haar measures on the cosets of $\mathrm{Ad}\, G$ in $\mathrm{Aut}(\mathrm{Ad}\, G)$.*

*Proof of Theorem* 4.1. Let $(\Omega, m)$ be a ME coupling of $\Gamma_1$ with $\Gamma_2$. Let $X \subset \Omega$ be a fundamental domain for $\Omega/\Gamma_2$, and

$$\alpha : \Gamma_1 \times X \to \Gamma_2$$

be the associated measurable cocycle. Recall, that by its definition

$$\gamma_1 \, x \in X \alpha(\gamma_1, x)$$

for a.e. $x \in X$ and all $\gamma_1 \in \Gamma_1$. Consider the cocycle

$$A : \Gamma_1 \times X \xrightarrow{\ \alpha\ } \Gamma_2 \xrightarrow{\ \mathrm{Ad}\ } \mathrm{Ad}\,\Gamma_2 \subset \mathrm{Ad}\,G$$

as an $\mathrm{Ad}\,G$-valued cocycle.

LEMMA 4.2. *The cocycle* $A : \Gamma_1 \times X \to \mathrm{Ad}\,G$ *is Zariski dense, i.e., $A$ is not measurably cohomologous to a cocycle taking values in a proper algebraic subgroup* $L \subset \mathrm{Ad}\,G$.

*Proof.* This is an adaptation of the argument due to Zimmer (c.f. [Zi3, p. 99]). Suppose that $A(\gamma, x) = \phi(\gamma \cdot x)^{-1} C(\gamma, x)\, \phi(x)$ for some measurable cocycle $C : \Gamma_1 \times X \to L$, where $L \subset \mathrm{Ad}\,G$ is a proper algebraic subgroup and $\phi : X \to \mathrm{Ad}\,G$ is a measurable map. Extend $\phi$ to the $\Gamma_2$-equivariant map $\Omega \to \mathrm{Ad}\,G$ by $\phi(x\gamma_2) = \phi(x)\,\mathrm{Ad}(\gamma_2)$ for $x \in X \subset \Omega$ and $\gamma_2 \in \Gamma_2$ (recall that $X \subset \Omega$ is a $\Gamma_2$ fundamental domain). By the definition of $\alpha$, for any $\gamma_1 \in \Gamma_1$ and a.e. $x \in X$:

$$\gamma_1 x = (\gamma_1 \cdot x)\, \alpha(\gamma_1, x),$$

so that using the extended definition of $\phi$:

$$
\begin{aligned}
\phi(\gamma_1 \, x) &= \phi((\gamma_1 \cdot x)\, \alpha(\gamma_1, x)) = \phi(\gamma_1 \cdot x)\, A(\gamma_1, x) \\
&= \phi(\gamma_1 \cdot x)\, \phi(\gamma_1 \cdot x)^{-1}\, C(\gamma_1, x)\, \phi(x) \\
&= C(\gamma_1, x)\, \phi(x).
\end{aligned}
$$

Thus, if $p : \mathrm{Ad}\,G \to L \backslash \mathrm{Ad}\,G$ is the projection, then the function

$$f : \Omega \xrightarrow{\ \phi\ } \mathrm{Ad}\,G \xrightarrow{\ p\ } L \backslash \mathrm{Ad}\,G$$

is $\Gamma_1$-invariant: $f(\gamma_1\omega) = f(\omega)$. Note also that $f$ is $\Gamma_2$-equivariant: $f(\omega\,\gamma_2) = f(\omega)\,\mathrm{Ad}(\gamma_2)$. Hence $f$ defines a $\Gamma_2$-equivariant map from $\Gamma_1 \backslash \Omega \cong Y$ to $L \backslash \mathrm{Ad}\,G$, which takes the finite measure $\nu$ to an $\mathrm{Ad}\,\Gamma_2$-invariant finite measure $\tilde{\nu}$ on $L \backslash \mathrm{Ad}\,G$. Setting

$$\bar{\nu} = \int_{\mathrm{Ad}\,\Gamma_2 \backslash \mathrm{Ad}\,G} \tilde{\nu}\, g \, dg$$

we obtain an $\mathrm{Ad}\,G$-invariant finite measure $\bar{\nu}$ on $L \backslash \mathrm{Ad}\,G$, which contradicts Borel's density theorem.                                              □

Let us assume, for a moment, that $\Gamma_1$ acts ergodically on $(X, \mu)$. Then applying the superrigidity for cocycles theorem ([Zi3, 5.2.5]) to the Zariski dense cocycle $A : \Gamma_1 \times X \to \mathrm{Ad}\, G$, we can conclude that there exists a homomorphism $\pi : \Gamma_1 \to \mathrm{Ad}\, G$ and a measurable map $\phi : X \to \mathrm{Ad}\, G$ so that

$$(4.2) \qquad A(\gamma, x) = \phi(\gamma \cdot x)^{-1} \pi(\gamma)\, \phi(x), \qquad (x \in X,\ \gamma \in \Gamma_1).$$

Since $A$ is Zariski dense, so is the image $\pi(\Gamma_1)$, and by Margulis's superrigidity [Ma] the homomorphism $\pi$ extends to an automorphism of $\mathrm{Ad}\, G$.

If $\pi$ is an inner automorphism, i.e. $\pi(h) = g^{-1}\, h\, g$, then replacing $\phi(x)$ by $g\, \phi(x)$, we can assume that $\pi(\gamma) = \mathrm{Ad}(\gamma)$ in (4.2). In general, $\pi$ is an inner automorphism in $\mathrm{Aut}(\mathrm{Ad}\, G)$, and replacing $\phi : X \to \mathrm{Ad}\, G$ by $\Phi = \pi\, \phi : X \to \mathrm{Aut}(\mathrm{Ad}\, G)$, we shall obtain

$$(4.3) \qquad A(\gamma, x) = \Phi(\gamma \cdot x)^{-1} \mathrm{Ad}(\gamma)\, \Phi(x)$$

in $\mathrm{Aut}(\mathrm{Ad}\, G)$. Coming back to the general measure-preserving (rather than ergodic) case, we can decompose $(\Omega, m)$ into ergodic components (Lemma 2.2) and proceed as above, so that $\Phi : X \to \mathrm{Aut}(\mathrm{Ad}\, G)$ will satisfy (4.3) for all $\gamma \in \Gamma_1$ and $m$-a.e. $x \in X$.

Recalling that $(X, \mu)$ is a $\Gamma_2$ fundamental domain, let us extend the map $\Phi : X \to \mathrm{Aut}(\mathrm{Ad}\, G)$ to $\Phi : \Omega \to \mathrm{Aut}(\mathrm{Ad}\, G)$ in a $\Gamma_2$-equivariant way:

$$\Phi(x\, \gamma_2) = \Phi(x)\, \mathrm{Ad}(\gamma_2), \qquad (x \in X \subset \Omega,\ \gamma_2 \in \Gamma_2).$$

*Claim* 4.3.  The map $\Phi : \Omega \to \mathrm{Aut}(\mathrm{Ad}\, G)$ is $\Gamma_1 \times \Gamma_2$-equivariant, i.e.

$$\Phi(\gamma_1\, \omega\, \gamma_2) = \mathrm{Ad}(\gamma_1)\, \Phi(\omega)\, \mathrm{Ad}(\gamma_2).$$

*Proof.*  By definition $\Phi$ is $\Gamma_2$-equivariant. It remains to show that

$$\Phi(\gamma_1\, \omega) = \mathrm{Ad}(\gamma_1)\, \Phi(\omega)$$

for $\gamma_1 \in \Gamma_1$. Let $\omega = x\, \gamma_2$ for some $x \in X \subset \Omega$ and $\gamma_2 \in \Gamma_2$. Then

$$\gamma_1\, \omega = \gamma_1\, x\, \gamma_2 = (\gamma_1 \cdot x)\, \alpha(\gamma_1, x)\, \gamma_2$$

with $\gamma_1 \cdot x \in X$ and $\alpha(\gamma_1, x)\, \gamma_2 \in \Gamma_2$. Therefore

$$
\begin{aligned}
\Phi(\gamma_1\, \omega) &= \Phi(\gamma_1 \cdot x)\, A(\gamma_1, x)\, \mathrm{Ad}(\gamma_2) \\
&= \Phi(\gamma_1 \cdot x)\, \Phi(\gamma_1 \cdot x)^{-1} \mathrm{Ad}(\gamma_1)\, \Phi(x)\, \mathrm{Ad}(\gamma_2) \\
&= \mathrm{Ad}(\gamma_1)\, \Phi(x)\, \mathrm{Ad}(\gamma_2) = \mathrm{Ad}(\gamma_1)\, \Phi(x\, \gamma_2) \\
&= \mathrm{Ad}(\gamma_1)\, \Phi(\omega).
\end{aligned}
$$

$\square$

This completes the proof of the existence of $\Gamma_1 \times \Gamma_2$-equivariant map $\Phi : \Omega \to \mathrm{Aut}(\mathrm{Ad}\, G)$.

PROPOSITION 4.4.  *The $\Gamma_1 \times \Gamma_2$-equivariant measurable map $\Phi : \Omega \to \mathrm{Aut}(\mathrm{Ad}\, G)$ is unique.*

*Remark* 4.5. Note that this, in particular, implies that the function $\phi : X \to \mathrm{Ad}\,G$, and the homomorphism $\pi$ are uniquely defined. It also follows that the definition of $\Phi : \Omega \to \mathrm{Aut}(\mathrm{Ad}\,G)$ does not depend on the choice of the fundamental domain $(X, \mu)$ from which $\alpha$, $\phi$ and $\pi$ were derived.

*Proof.* Suppose $\Phi$ and $\Phi'$ are two measurable $\Gamma_1 \times \Gamma_2$-equivariant maps. Define $\Psi : \Omega \to \mathrm{Aut}(\mathrm{Ad}\,G)$ by $\Psi(\omega) = \Phi'(\omega)\,\Phi(\omega)^{-1}$. Observe that $\Psi$ satisfies

$$\Psi(\gamma_1\,\omega\,\gamma_2) = \mathrm{Ad}(\gamma_1)\,\Psi(\omega)\,\mathrm{Ad}(\gamma_1)^{-1} \qquad (\gamma_1 \in \Gamma_1,\ \gamma_2 \in \Gamma_2).$$

In particular, $\Psi$ is right $\Gamma_2$-invariant; hence it can be considered as a function on a $\Gamma_2$-fundamental domain $(X, \mu)$. With respect to the $\Gamma_1$-action on $(X, \mu)$ $(\gamma_1 : x \mapsto \gamma_1 \cdot x)$ we have

$$\Psi(\gamma \cdot x) = \mathrm{Ad}(\gamma)\,\Psi(x)\,\mathrm{Ad}(\gamma)^{-1}, \qquad (\gamma \in \Gamma_1).$$

We claim that $m$-a.e. $\Psi(\omega) = e$. This can be deduced from a more general Lemma 5.3. However, in the present case the idea of the general argument can be presented in an almost elementary way as follows.

It is enough to show that for any compact $K \subset \mathrm{Aut}(\mathrm{Ad}\,G)$ we have $\Psi(x) = e$ for $\mu$-a.e. $x \in E_K = \Psi^{-1}(K)$. Fix $\gamma \in \Gamma$. By Poincaré recurrence theorem $\mu$-a.e. $x \in E_K$ returns to $E_K$ infinitely often, so that

$$\Psi(\gamma^n \cdot x) = \mathrm{Ad}(\gamma)^n\,\Psi(x)\,\mathrm{Ad}(\gamma)^{-n} \in K$$

infinitely often. It is well known that for any fixed regular element $g \in \mathrm{Ad}\,G$ and any $h \in \mathrm{Ad}\,G$ which does not commute with $g$, one has $g^n\,h\,g^{-n} \to \infty$ as $n \to \infty$ or $n \to -\infty$. The same holds for $h \in \mathrm{Aut}(\mathrm{Ad}\,G)$. Hence, for $\mu$-a.e. $x \in E_K$ the element $\Psi(x)$ commutes with all $\gamma \in \mathrm{Ad}\,\Gamma_1$. Borel's density theorem implies that on $E_K$ a.e. $\Psi(x)$ belongs to the center of $\mathrm{Aut}(\mathrm{Ad}\,G)$ which is trivial. $\qquad\qquad\square$

Now let $D \subset \mathrm{Aut}(\mathrm{Ad}\,G)$ be a $\Gamma_2$-fundamental domain. Since $\Phi$ is $\Gamma_1 \times \Gamma_2$-equivariant, the set $X = \Phi^{-1}(D)$ is a $\Gamma_2$-fundamental domain in $\Omega$. Thus the projection $\Phi_*\mu$ of the finite measure $\mu = m|_X$ to $D$ is an $\mathrm{Ad}\,\Gamma_1$-invariant finite measure on $D \cong \mathrm{Aut}(\mathrm{Ad}\,G)/\mathrm{Ad}\,\Gamma_2$ and, by disintegrating $\mu$ (and $m$) with respect to $\Phi_*\mu$ (and $\Phi_*m$) one obtains probability measures on the fibers.

$\mathrm{Aut}(\mathrm{Ad}\,G)$ consists of a finite number of $\mathrm{Ad}\,G$-cosets, each of which is $\mathrm{Ad}\,\Gamma_1 \times \mathrm{Ad}\,\Gamma_2$-invariant. Restricting $\Phi_*m$ to each of $\mathrm{Ad}\,G$-cosets one easily deduces the last statement of Theorem 4.1 from the following:

LEMMA 4.6. *Let $G$ be a simple connected Lie group with trivial center, and let $\Gamma_1$, $\Gamma_2 \subset G$ be lattices. Suppose that $\mu$ is a probability Borel measure on $G/\Gamma_2$, which is invariant and ergodic under the left $\Gamma_1$-action on $G/\Gamma_2$. Then either $\mu$ is a finite atomic measure, or $\mu$ is the $G$-invariant probability measure on $G/\Gamma_2$.*

*Proof.* The assertion follows from Ratner's theorem. Let $\tilde{\mu}$ be the (right) $\Gamma_2$-invariant lifting of $\mu$ from $G/\Gamma_2$ to $G$. Let $m = m_G$ denote a bi-invariant Haar measure on $G$. Consider the measure $\tilde{M}$ on $G \times G$, defined by

$$\int_{G \times G} f(g_1, g_2) \, d\tilde{M}(g_1, g_2) = \int_G \int_G f(g_1, g_2) \, d\tilde{\mu}(g_1^{-1} g_2) \, dm(g_1).$$

Note that $\tilde{M}$ is invariant under the (right) actions of $\Gamma_1 \subset G \times \{e\}$ and $\Gamma_2 \subset \{e\} \times G$, so $\tilde{M}$ is a lifting of a measure $M$ on $G/\Gamma_1 \times G/\Gamma_2$. It is easily seen that

$$M = \int_{G/\Gamma_1} g^{-1} \mu \, dm(g)$$

is a finite measure. Moreover, $\tilde{M}$, and hence $M$, are $\Delta(G)$-invariant, where

$$\Delta(G) := \{(g, g) \in G \times G \mid g \in G\}.$$

Since $\Gamma_1 \times \Gamma_2$ forms a lattice in $G \times G$, and $\Delta(G) \subset G \times G$ is generated by unipotents, Ratner's theorem [Ra] implies that $M$ is supported on an orbit of a closed subgroup $L \subset G \times G$, where the intersection $\Gamma_L = (\Gamma_1 \times \Gamma_2) \cap L$ is a lattice in $L$ and $\tilde{M} = m_L$ is a Haar measure on $L$. Let $L_0 \subseteq L$ be the connected component of the identity.

LEMMA 4.7. *Let $G$ be a simple connected Lie group with finite center. Let $L_0$ be a connected subgroup of $G \times G$ which contains the diagonal $\Delta(G) = \{(g, g) \in G \times G \mid g \in G\}$. Then either $L_0 = \Delta(G)$ or $L_0 = G \times G$.*

*Proof.* For $g \in G$, let $L(g) \subseteq G$ be the fiber of $L_0$ over $g$, i.e. $L(g) \times \{g\} = L_0 \cap G \times \{g\}$. Then $L(e)$ is a closed subgroup of $G$. Note that for any $g \in G$:

$$\Delta(g) \cdot (L_0 \cap G \times \{e\}) = L_0 \cap G \times \{g\} = (L_0 \cap G \times \{e\}) \cdot \Delta(g).$$

Hence $g \, L(e) = L(g) = L(e) \, g$, so that $L(e)$ is a normal closed subgroup of $G$. Therefore, either $L(e) = \{e\}$ and $L(g) = \{g\}$, or $L(e) = L(g) = G$. In the former case $L_0 = \Delta(G)$ and in the latter case $L_0 = G \times G$. $\qquad \square$

It is easily verified that in the case of $L_0 = \Delta(G)$ the measure $\mu$ is finite atomic; while in the case of $L_0 = G \times G$ the measure $\mu$ is the unique $G$-invariant probability measure on $G/\Gamma_2$. This proves Lemma 4.6. $\qquad \square$

## 5. Construction of the representation $\rho : \Lambda \to \mathrm{Aut}(\mathrm{Ad}\, G)$

This crucial step of the proof of Theorem 3.1 describes the construction of a family of mutually equivalent representations $\rho_x : \Lambda \to \mathrm{Aut}(\mathrm{Ad}\, G)$ of an unknown group $\Lambda$ which is ME to a lattice $\Gamma$. Let $(\Sigma, \sigma)$ be a ME coupling of $\Gamma$ with $\Lambda$. The idea is to use the map $\Phi_\Omega$ constructed in Theorem 4.1 from the $\Gamma$-self ME coupling $\Omega = (\Sigma \times_\Lambda \Lambda \times_\Lambda \tilde{\Sigma})$ to $\mathrm{Aut}(\mathrm{Ad}\, G)$.

*Example* 5.1.   Consider the case where both $\Gamma$ and $\Lambda$ are lattices in $G$ with $(G, m_G)$ being the coupling. Assume that $G$ has trivial center. Then the map

$$\Phi : (G \times_\Lambda \Lambda \times_\Lambda \check{G}) \to G \cong \operatorname{Ad} G \subset \operatorname{Aut}(\operatorname{Ad} G)$$

which is given by $\Phi([x, \lambda, \check{y}]) = x \, \lambda \, y^{-1}$ is the one constructed in Theorem 4.1. Observe that in this case the map

$$\Phi([x, \lambda, \check{y}]) \, \Phi([x, e, \check{y}])^{-1} = (x \, \lambda \, y^{-1}) \, (x \, y^{-1})^{-1} = x \, \lambda \, x^{-1}$$

does not depend on $y$, and for a.e. fixed $x$ defines a representation of $\Lambda$.

*Preliminaries.* We shall prove that somewhat similar phenomenon exists in the general case of an unknown $\Lambda$. The following lemma describes how one constructs a representation, given a measurable function satisfying certain a.e. identities. Note that the data consists of some identities on $\Sigma^n$ which hold $\sigma^n$-almost everywhere, but it is not known whether these identities hold *everywhere* on $S^n$ for any $\sigma$-conull $S \subseteq \Sigma$. This is a common feature of the measure-theoretic framework.

Lemma 5.2.   *Let a countable group $\Lambda$ act (from the right) on a measure space $(\Sigma, \sigma)$, and let $F : \Sigma \times \Sigma \to G$ be some measurable map to a* lcsc *group $G$. If $F$ satisfies*

$(C_{\mathrm{inv}})$    $F(x, y) = F(x\lambda, \, y\lambda)$    $\sigma^2$-*a.e. on* $\Sigma^2$ *for all* $\lambda \in \Lambda$.

$(C_{\mathrm{cncl}})$    $F(x\lambda, \, y) \, F(x, y)^{-1} = F(x\lambda, \, z) \, F(x, z)^{-1}$    $\sigma^3$-*a.e. on* $\Sigma^3$.

*Then a.e. $x \in \Sigma$ defines a homomorphism $\rho_x : \Lambda \to G$, given by*

$$\rho_x(\lambda) = F(x\lambda, y) \, F(x, y)^{-1}, \qquad (\lambda \in \Lambda).$$

*If, moreover, $F(x, y)$ satisfies*

$(C_{\mathrm{sym}})$    $F(x, y) = F(y, x)^{-1}$    $\sigma^2$-*a.e. on* $\Sigma^2$

$(C_{\mathrm{coc}})$    $F(x, y) \, F(y, z) = F(x, z)$    $\sigma^3$-*a.e. on* $\Sigma^3$

*then for $\sigma^2$-a.e. $(x, y) \in \Sigma^2$ the representations $\rho_x$ and $\rho_y$ are equivalent:*

$$\rho_y(\lambda) = F(x, y)^{-1} \, \rho_x(\lambda) \, F(x, y), \qquad (\lambda \in \Lambda).$$

*Proof.* By $C_{\mathrm{cncl}}$ for $\sigma$-a.e. $x \in \Sigma$ the function $\rho_x(\lambda) = F(x\lambda, y) \, F(x, y)^{-1}$ has the same value for a.e. $y \in \Sigma$. For any $\lambda_1, \lambda_2 \in \Lambda$ and a.e. $x \in \Sigma$, choosing a.e. $y \in \Sigma$ and using $C_{\mathrm{inv}}$ and $C_{\mathrm{cncl}}$, we have:

$$
\begin{aligned}
\rho_x(\lambda_1 \, \lambda_2^{-1}) &= F(x\lambda_1\lambda_2^{-1}, \, y) \, F(x, y)^{-1} \\
&= F(x\lambda_1, \, y\lambda_2) \, F(x\lambda_2, \, y\lambda_2)^{-1} \\
&= F(x\lambda_1, y\lambda_2) \, F(x, \, y\lambda_2)^{-1} \left( F(x\lambda_2, \, y\lambda_2) \, F(x, \, y\lambda_2)^{-1} \right)^{-1} \\
&= \rho_x(\lambda_1) \, \rho_x(\lambda_2)^{-1}.
\end{aligned}
$$

Since $\Lambda$ is countable, we conclude that for $\sigma$-a.e. $x$, the map $\rho_x : \Lambda \to G$ is a homomorphism. Assume that $F(x, y)$ satisfies also the conditions $C_{\mathrm{coc}}$ and $C_{\mathrm{sym}}$. Then for a.e. $(x, y) \in \Sigma^2$ we have for a.e. $z \in \Sigma$:

$$
\begin{aligned}
F(x, y)^{-1} \rho_x(\lambda) F(x, y) &= F(x, y)^{-1} F(x\lambda, z) F(x, z)^{-1} F(x, y) \\
&= F(x\lambda, y\lambda)^{-1} F(x\lambda, z) F(x, z)^{-1} F(x, y) \\
&= F(y\lambda, z) F(y, z)^{-1} = \rho_y(\lambda). \qquad \square
\end{aligned}
$$

The crucial condition to be verified for an application of Lemma 5.2 is $C_{\mathrm{cncl}}$. The proof of this property for an appropriately chosen function $F(x, y)$ will rely on the following lemmas:

LEMMA 5.3. *Let a group $\Lambda$ act ergodically (from the right) on a finite measure space $(W, \eta)$, let $G$ be a semisimple Lie group with trivial center, and $B : W \times \Lambda \to G$ be a measurable cocycles, which is Zariski dense in $G$. Suppose, moreover, that there is a measurable map $M : W \to \mathcal{P}(G)$ with values in the space of probability measures $\mathcal{P}(G)$ on $G$, which satisfies*

$$
M(w \cdot \lambda) = B(w, \lambda)^{-1} M(w) B(w, \lambda), \quad (\lambda \in \Lambda)
$$

*then for $\eta$-a.e. $M(w) = \delta_e$ is a the point measure at the identity $e \in G$.*

*Proof.* The group $G$ acts on itself by conjugation: $\mathrm{Ad}(g) : g' \mapsto g^{-1} g' g$. Since the center is trivial, $\mathrm{Ad} : G \to \mathrm{Ad}\, G$ is an isomorphism. Moreover since the $\mathrm{Ad}\, G$-action on itself is essentially algebraic, the corresponding action on the space of probability measures $\mathcal{P}(G)$ is smooth (see [Zi3, 3.2.6]), in the sense that there exists a countable family of $\mathrm{Ad}\, G$-invariant measurable sets in $\mathcal{P}(G)$ which separate orbits. Without loss of generality we can assume that $\Lambda$ acts ergodically on $W$. Then Zimmer's cocycle reduction lemma (see [Zi3]) implies that $M$ is supported on a single $\mathrm{Ad}\, G$-orbit: $M(w) = \mathrm{Ad}(f(w))M_0$, and $B(g, x)$ is cohomologous to a cocycle $B_0$, taking values in a stabilizer of $M_0 \in \mathcal{P}(G)$. Since the stabilizers are algebraic ([Zi3, 3.2.4]) and $B(g, x)$ is assumed to be Zariski dense, $M_0$ is $\mathrm{Ad}\, G$-invariant. This implies that $M_0 = \delta_e$, and thus $M(w) = \delta_e$ a.e. on $W$. $\qquad \square$

LEMMA 5.4. *Let $\Gamma$ be a lattice in a semisimple Lie group $G$. Suppose $\Lambda$ is some group which is ME to $\Gamma$, let $(\Sigma, \sigma)$ be their ME coupling, let $(X, \mu)$, $(Y, \nu)$ be $\Lambda$- and $\Gamma$-fundamental domains, and let $\alpha : \Gamma \times X \to \Lambda$ and $\beta : Y \times \Lambda \to \Gamma$ be the associated cocycles. Then the cocycle*

$$
B : Y \times \Lambda \xrightarrow{\beta} \Gamma \xrightarrow{\mathrm{Ad}} \mathrm{Ad}\, \Gamma \subset \mathrm{Ad}\, G
$$

*is Zariski dense in $\mathrm{Ad}\, G$.*

*Proof.* The proof follows from Lemma 4.2, applied to the cocycle associated to the composition coupling $\Sigma \times_\Lambda \tilde{\Sigma}$ of $\Gamma$ with $\Gamma$. $\qquad \square$

*First candidate for $F(x, y)$.* Let $(\Sigma, \sigma)$ be a ME coupling between two countable groups $\Gamma$ and $\Lambda$, where $\Gamma$ is a lattice in $G$ as in Theorem 3.1. Consider the dual coupling $(\check{\Sigma}, \check{\sigma})$ of $\Lambda$ with $\Gamma$, and let $(\Omega, m)$ be the composition coupling of $\Gamma$ with $\Gamma$:

$$(\Omega, m) = (\Sigma \times_\Lambda \Lambda \times_\Lambda \check{\Sigma}, \ \sigma \times_\Lambda m_\Lambda \times_\Lambda \check{\sigma}).$$

Let $\Phi_\Omega : \Omega = (\Sigma \times_\Lambda \Lambda \times_\Lambda \check{\Sigma}) \to \operatorname{Aut}(\operatorname{Ad} G)$ be the measurable map as in Theorem 4.1. Let us simplify the notations introducing a measurable map $\tilde{\Phi} : \Sigma \times \Sigma \to \operatorname{Aut}(\operatorname{Ad} G)$, defined by

$$(5.1) \qquad\qquad\qquad \tilde{\Phi}(x, y) = \Phi_\Omega([x, e, \check{y}]).$$

By its definition $\tilde{\Phi}$ satisfies

$$(5.2) \qquad \tilde{\Phi}(x\lambda, \ y\lambda) = \Phi_\Omega([x\lambda, e, \lambda^{-1}\check{y}]) = \Phi_\Omega([x, \lambda\lambda^{-1}, \check{y}]) = \tilde{\Phi}(x, \ y)$$

while $\Gamma \times \Gamma$-equivariance of $\Phi_\Omega$ gives

$$(5.3) \qquad \begin{aligned} \tilde{\Phi}(\gamma x, \ y) &= \operatorname{Ad}(\gamma)\,\tilde{\Phi}(x, \ y), && (\gamma \in \Gamma) \\ \tilde{\Phi}(x, \ \gamma y) &= \tilde{\Phi}(x, \ y)\operatorname{Ad}(\gamma)^{-1}, && (\gamma \in \Gamma). \end{aligned}$$

Having Example 5.1 and Lemma 5.2 in mind, an optimistic reader would expect function $\tilde{\Phi}(x, y)$ to satisfy the conditions of Lemma 5.2. The crucial property to be verified is $C_{\mathrm{cncl}}$, i.e. that

$$\Phi([x, \lambda, \check{y}])\,\Phi([x, e, \check{y}])^{-1} = \tilde{\Phi}(x\lambda, \ y)\,\tilde{\Phi}(x, y)^{-1}$$

does not depend on $y$. Unfortunately, we cannot show this directly, although by the end of the proof we shall see (Remark 7.2) that $\tilde{\Phi}(x, y) = \Phi_\Sigma(x)\,\Phi_\Sigma(y)^{-1}$ for some measurable $\Phi_\Sigma : \Sigma \to \operatorname{Aut}(\operatorname{Ad} G)$, so that $\tilde{\Phi}$ indeed satisfies all the conditions in Lemma 5.2.

*The choice of $F(x, y)$ which works.* At this point we choose to consider another map $\Psi : \Sigma \times \Sigma \times \Sigma \to \operatorname{Aut}(\operatorname{Ad} G)$ defined by

$$(5.4) \qquad\qquad\qquad \Psi(u, x, y) = \tilde{\Phi}(u, y)\,\tilde{\Phi}(x, y)^{-1}.$$

We shall show that $\Psi(u, x, y)$ does not depend on $y$, namely $\Psi(u, x, y) = \Psi(u, x, z)$, $\sigma^4$-a.e. on $\Sigma^4$. (Note, that this still does not prove $C_{\mathrm{cncl}}$ for $\tilde{\Phi}$, because the latter is given in terms of a zero measure set in $\Sigma^4$.)

Properties (5.2) and (5.3) imply that $\Psi$ satisfies:

$$(5.5) \qquad \begin{aligned} \Psi(\gamma x, \ y, \ z) &= \operatorname{Ad}(\gamma)\,\Psi(x, \ y, \ z) && (\gamma \in \Gamma) \\ \Psi(x, \ \gamma y, \ z) &= \Psi(x, \ y, \ z)\operatorname{Ad}(\gamma)^{-1} && (\gamma \in \Gamma) \\ \Psi(x, \ y, \ \gamma z) &= \Psi(x, \ y, \ z) && (\gamma \in \Gamma) \\ \Psi(x\lambda, \ y\lambda, \ z\lambda) &= \Psi(x, \ y, \ z) && (\lambda \in \Lambda). \end{aligned}$$

The following lemma is the key point of the construction of the representation:

LEMMA 5.5. *The map* $\Psi : \Sigma \times \Sigma \times \Sigma \to G$, *defined by (5.4) does not depend on the third coordinate, i.e.* $\Psi(x, y, z_1) = \Psi(x, y, z_2)$ $\sigma^4$-*a.e. on* $\Sigma^4$.

*Proof.* Let $Y \subset \Sigma$ be a fundamental domain for the left $\Gamma$-action on $\Sigma$, and let $\beta : \Sigma \times \Lambda \to \Gamma$ be the associated cocycle. The right action of $\Lambda$ on $\Gamma \backslash \Sigma \cong Y$ is given by

$$y \cdot \lambda = \beta(y, \lambda)^{-1} y \lambda, \qquad (y \in Y, \ \lambda \in \Lambda)$$

where the left $\Gamma$-action and the right $\Lambda$-actions on the right-hand side are in $\Sigma$. By (5.5), it is enough to show that $\Psi$, restricted to $Y \times Y \times Y$, does not depend on the third coordinate. Denote by $B = \mathrm{Ad} \circ \beta$ the cocycle

$$B : Y \times \Lambda \xrightarrow{\ \beta\ } G \xrightarrow{\ \mathrm{Ad}\ } \mathrm{Ad}\, G \subset \mathrm{Aut}(\mathrm{Ad}\, G).$$

Identities (5.5) yield the following crucial relation:

$$(5.6) \quad \begin{aligned} \Psi(x \cdot \lambda, \ y \cdot \lambda, \ z \cdot \lambda) &= \Psi\left(\beta(x, \lambda)^{-1} x \lambda, \ \beta(y, \lambda)^{-1} y \lambda, \ \beta(z, \lambda)^{-1} z \lambda\right) \\ &= B(x, \lambda)^{-1} \Psi(x\lambda, \ y\lambda, \ z\lambda) \, B(y, \lambda) \\ &= B(x, \lambda)^{-1} \Psi(x, \ y, \ z) \, B(y, \lambda). \end{aligned}$$

We shall now use the fact that the transformation: $\Psi \mapsto B(x, \lambda)^{-1} \Psi B(y, \lambda)$ in (5.6) does not involve the $z$-variable. Let $M(x, y) \in \mathcal{P}(\mathrm{Aut}(\mathrm{Ad}\, G))$ be the distribution of $\Psi(x, y, z_1) \Psi(x, y, z_2)^{-1}$ as $z_1, z_2 \in Y$ are chosen independently with $\nu$-distribution. In other words, for $x, y \in Y$ define $M(x, y)$ to be

$$dM(x, y) = \Psi(x, y, z_1) \Psi(x, y, z_2)^{-1} \, d\nu(z_1) \, d\nu(z_2).$$

By (5.6) we have

$$\begin{aligned} \Psi(x \cdot \lambda, \ y \cdot \lambda, \ z_1 \cdot \lambda) \, \Psi(x \cdot \lambda, \ y \cdot \lambda, \ z_2 \cdot \lambda)^{-1} \\ = B(x, \lambda)^{-1} \Psi(x, y, z_1) \, \Psi(x, y, z_2)^{-1} B(x, \lambda) \end{aligned}$$

and therefore, using $\Lambda$-invariance of $\nu$,

$$M(x \cdot \lambda, \ y \cdot \lambda) = B(x, \lambda)^{-1} M(x, y) B(x, \lambda).$$

Conjugation by $B(x, \lambda) \in \mathrm{Ad}\, G$ preserves the (finite number of) $\mathrm{Ad}\, G$ cosets in $\mathrm{Aut}(\mathrm{Ad}\, G)$. Thus we can assume that $M(x, y)$ is supported on $\mathrm{Ad}\, G$.

Now, taking $(W, \eta) = (Y \times Y, \nu \times \nu)$ with the diagonal $\Lambda$-action and viewing $B(x, \lambda)$ as a cocycle on $(W, \eta)$, we observe, that by Lemma 5.4, the cocycle $B(w, \lambda)$ is Zariski dense in the center free simple group $\mathrm{Ad}\, G$. Thus Lemma 5.3 shows that $M(x, y) = \delta_e$. $\qquad\square$

*Claim* 5.6. The function $F(x, y) = \Psi(x, y, z) = \tilde{\Phi}(x, z) \tilde{\Phi}(y, z)^{-1}$ satisfies the conditions $C_{\mathrm{inv}}$, $C_{\mathrm{cncl}}$, $C_{\mathrm{sym}}$ and $C_{\mathrm{coc}}$ of Lemma 5.2.

*Proof.* Condition $C_{inv}$ follows from $\Lambda$-invariance of $\tilde{\Phi}$. To verify property $C_{cncl}$, note that $\sigma^3$-a.e. on $\Sigma^3$

$$\begin{aligned}
F(x\lambda, y)\, F(x,y)^{-1} &= \tilde{\Phi}(x\lambda, z)\, \tilde{\Phi}(y,z)^{-1} \left( \tilde{\Phi}(x,z)\, \tilde{\Phi}(y,z)^{-1} \right)^{-1} \\
&= \tilde{\Phi}(x\lambda,\, z)\, \tilde{\Phi}(x,\, z)^{-1}.
\end{aligned}$$

The right-hand side does not depend on $y$, while the left-hand side does not depend on $z$. Hence $C_{cncl}$ follows. By its definition $F(x,y)$ satisfied $C_{sym}$. Now observe that $\sigma^4$-a.e. on $\Sigma^4$ we have

$$\begin{aligned}
F(x,y)\, F(y,z) &= \tilde{\Phi}(x,w)\, \tilde{\Phi}(y,w)^{-1}\, \tilde{\Phi}(y,w)\, \tilde{\Phi}(z,w)^{-1} \\
&= \tilde{\Phi}(x,w)\, \tilde{\Phi}(z,w)^{-1} = F(x,z),
\end{aligned}$$

which verifies condition $C_{coc}$.                                              $\square$

We can now apply Lemma 5.2 to produce a family of mutually equivalent homomorphisms $\rho_x : \Lambda \to G$.

## 6. The kernel of $\rho : \Lambda \to G$ is finite

LEMMA 6.1. *The subgroup* $\Lambda_0 = \mathrm{Ker}(\rho_x : \Lambda \to \mathrm{Aut}(\mathrm{Ad}\, G))$ *is at most finite.*

*Proof.* First note that since a.e. $\rho_x$ are equivalent, the group $\Lambda_0 = \mathrm{Ker}(\rho_x)$ is a well-defined normal subgroup in $\Lambda$, which does not depend on $x$. Hence the definition of $\rho_x$ yields that $\lambda \in \Lambda_0$ if and only if $F(x\lambda, y) = F(x,y)$ for a.e. $(x,y) \in \Sigma \times \Sigma$. Since for a.e. $z$, $F(x\lambda, y) = \tilde{\Phi}(x\lambda, z)\, \tilde{\Phi}(y,z)^{-1}$ and $F(x,y) = \tilde{\Phi}(x,z)\, \tilde{\Phi}(y,z)^{-1}$ simultaneously, we obtain

$$(6.1) \qquad \lambda \in \Lambda_0 \qquad \text{if and only if} \qquad \Phi_\Omega([x, \lambda, \check{z}]) = \Phi_\Omega([x, e, \check{z}]).$$

Let $D \subset \mathrm{Aut}(\mathrm{Ad}\, G)$ be a fundamental domain for $\mathrm{Aut}(\mathrm{Ad}\, G)/\mathrm{Ad}\,\Gamma$ and let $E = \Phi_\Omega^{-1}(D) \subset \Omega$ be its preimage. $\Gamma \times \Gamma$-equivariance of $\Phi_\Omega$ implies that $E$ forms a fundamental domain for $\Omega/\Gamma$.

Now observe that $\Omega = \Sigma \times_\Lambda \Lambda \times_\Lambda \check{\Sigma}$ can be represented in the form

$$\Omega = X \times \Lambda \times \check{X}$$

where $X$ is some fundamental domain for $\Sigma/\Lambda$. So $E$ is a disjoint union of $E(\lambda) = E \cap (X \times \{\lambda\} \times \check{X})$. The relation (6.1) implies that $E(\lambda) = E(\lambda\lambda_0)$ for any $\lambda_0 \in \Lambda_0$. Since $E = \bigcup_{\lambda \in \Lambda} E(\lambda)$ is a countable union, which has a finite positive measure, we deduce that $\Lambda_0$ has to be finite.        $\square$

## 7. End of the proof of Theorem 3.1

At this point we have proved that if $\Gamma \subset G$ is a lattice as in Theorem 3.1, and $\Lambda$ is a group which is ME to $\Gamma$, then there exists a representation $\rho : \Lambda \to \operatorname{Aut}(\operatorname{Ad} G)$ with finite kernel $\Lambda_0 = \operatorname{Ker}(\rho)$. Restricting $\rho$ to a finite index subgroup

$$\Lambda' = \operatorname{Ker}(\Lambda \xrightarrow{\rho} \operatorname{Aut}(\operatorname{Ad} G) \xrightarrow{o} \operatorname{Out}(\operatorname{Ad} G)),$$

we can assume that $\Lambda_2 = \rho(\Lambda') \subseteq \operatorname{Ad} G$. This enables us to apply the result of Zimmer ([Zi2, Cor. 1.2]) which shows that in this case, $\Lambda_2$ is a lattice in $\operatorname{Ad} G$. In fact, with our present setup, we need only Lemmas 2.6 and 2.7 of [Zi2] to be applied to the cocycle

$$\rho \circ \alpha : \ \Gamma \times \Sigma/\Lambda' \to \Lambda' \to \operatorname{Ad} G.$$

The first one ([Zi2, 2.6]), which is based on superrigidity for cocycles theorem (both real and $p$-adic cases), shows that $\Lambda_2 = \rho(\Lambda') \subset \operatorname{Ad} G$ is discrete.

The second one ([Zi2, 2.7]) shows that $\operatorname{Ad} G/\Lambda_2$ carries a finite invariant measure. Hence $\Lambda_2 = \rho(\Lambda')$ is a lattice in $\operatorname{Ad} G$, and therefore $\Lambda_1 = \rho(\Lambda)$ is a lattice in $\operatorname{Aut}(\operatorname{Ad} G)$. This proves the first part of Theorem 3.1.

Now consider $\Gamma_2 = \operatorname{Ad} \Gamma$ and $\Lambda_2 = \rho(\Lambda')$ as lattices in $\operatorname{Ad} G$. Let $(\Sigma_2, \sigma_2)$ be the factor space of $(\Sigma, \sigma)$ divided by the action of the finite group $(\Gamma_2 \cap Z(G)) \times (\Lambda_0 \cap \Lambda')$. It forms a ME coupling of $\Gamma_2$ and $\Lambda_2$-lattices in $\operatorname{Ad} G$. Applying Theorem 4.1 we obtain a (unique) measurable $\Gamma_2 \times \Lambda_2$-equivariant map $: \Sigma_2 \to \operatorname{Aut}(\operatorname{Ad} G)$. Lifting it to $\Sigma$ we obtain a measurable map $\Phi_\Sigma : \Sigma \to \operatorname{Aut}(\operatorname{Ad} G)$ satisfying $\sigma$-a.e.

$$(7.1) \qquad \Phi_\Sigma(\gamma\, x\, \lambda') = \operatorname{Ad}(\gamma)\, \Phi_\Sigma(x)\, \rho(\lambda') \qquad (\gamma \in \Gamma, \ \lambda' \in \Lambda').$$

*Claim* 7.1. The map $\Phi_\Sigma : \Sigma \to \operatorname{Aut}(\operatorname{Ad} G)$ is $\Gamma \times \Lambda$-equivariant (and not just $\Gamma \times \Lambda'$-equivariant).

*Proof.* Consider a function $f : \Sigma \times \Lambda \to \operatorname{Aut}(\operatorname{Ad} G)$ given by

$$f(x, \lambda) = \Phi_\Sigma(x\, \lambda)\, \left(\Phi_\Sigma(x)\, \rho(\lambda)\right)^{-1}.$$

Then for any $\lambda \in \Lambda$ and any $\lambda' \in \Lambda'$:

$$f(x\, \lambda',\, \lambda) = f(x,\, \lambda'\lambda) = f(x,\, \lambda\lambda').$$

So $f$ is actually defined on $X' \times \Lambda/\Lambda'$, where $X'$ is a $\Lambda'$-fundamental domain. If $M(x)$ is the uniform distribution of $f(x, \lambda)$ over $\Lambda/\Lambda'$, then

$$M(\gamma \cdot x) = \operatorname{Ad}(\gamma)\, M(x)\, \operatorname{Ad}(\gamma)^{-1} \qquad (\gamma \in \Gamma)$$

and Lemma 5.3 (or the similar argument in Proposition 4.4) implies that $M(x) = e$, so that $\Phi_\Sigma(x\lambda) = \Phi(x)\, \rho(\lambda)$ for all $\lambda \in \Lambda$. $\square$

This completes the proof of Theorem 3.1.

*Remark* 7.2. In a retrospective on the proof, one can see (using the uniqueness part of Theorem 4.1) that the map $\Phi_\Omega : (\Sigma \times_\Lambda \Lambda \times_\Lambda \tilde{\Sigma}) \to \mathrm{Aut}(\mathrm{Ad}\,G)$ in the proof of Theorem 3.1 is given by

$$\Phi_\Omega([x, \lambda, \check{y}]) = \Phi_\Sigma(x)\,\rho(\lambda)\,\Phi_\Sigma(y)^{-1}$$

and therefore the maps $\tilde{\Phi}$, $F : \Sigma \times \Sigma \to G$ in (5.1) and in Claim 5.6 satisfy

$$F(x, y) = \tilde{\Phi}(x, y) = \Phi_\Sigma(x)\,\Phi_\Sigma(y)^{-1}$$

and finally,

$$\rho_x(\lambda) = \Phi_\Sigma(x)\,\rho(\lambda)\,\Phi_\Sigma(x)^{-1}.$$

## 8. Measure Equivalence and unitary representations

Let $(\Omega, m)$ be a ME coupling of two countable groups $\Gamma$, $\Lambda$ and let $(X, \mu)$ and $(Y, \nu)$ be $\Lambda$- and $\Gamma$-fundamental domains in $\Omega$. Denote by $\tilde{V}$ the $L^2$ space of measurable functions $X \to V$. More precisely,

$$\tilde{V} = L^2(X, V) = \left\{ f : X \to V \,\Big|\, \int_X \|f(x)\|^2\,d\mu(x) < \infty \right\}.$$

Let $\Gamma$ act on $\tilde{V}$ by

$$(\tilde{\pi}(\gamma)f)(x) = \pi(\alpha(\gamma^{-1}, x)^{-1})\left(f(\gamma^{-1} \cdot x)\right).$$

One checks that $(\tilde{\pi}, \tilde{V})$ is a unitary $\Gamma$-representation, which will be called the *induced* representation.

Recall that a unitary representation $(\pi, V)$ of a discrete group $\Lambda$ is said to contain *almost invariant* vectors, if there exists a sequence $\{v_n\}$ of unit vectors such that $\|\pi(\lambda)v_n - v_n\| \to 0$ as $n \to \infty$ for every $\lambda \in \Lambda$.

LEMMA 8.1. *If the $\Lambda$-representation $(\pi, V)$ contains almost invariant vectors, then so does the induced $\Gamma$-representation $(\tilde{\pi}, \tilde{V})$.*

*Proof.* Let $\{v_n\}$ be a sequence in $V$ of almost $\Lambda$-invariant unit vectors and let $f_n \in \tilde{V}$ be the sequence of vectors $f_n(x) \equiv v_n$. Then $\|f_n\| = \mu(X)^{1/2}$, and for any fixed $\gamma \in \Gamma$

$$\langle \tilde{\pi}(\gamma)\,f_n,\,f_n \rangle = \int_X \left\langle \alpha(\gamma^{-1}, x)^{-1}\,v_n,\,v_n \right\rangle d\mu(x).$$

For a sufficiently large finite set $F \subset \Lambda$, one has $\alpha(\gamma^{-1}, x)^{-1} \in F$ on $X \setminus E$ with $\mu(E)$ being small. As $n \to \infty$, the vectors $v_n$ become almost $F$-invariant, so the integrand is close to one on $X \setminus E$, and is bounded by 1 on $E$. This shows that $(\tilde{\pi}, \tilde{V})$ contains $\Gamma$-almost invariant vectors. $\qquad\square$

A unitary $\Lambda$-representation $(\pi, V)$ is said to be *weakly mixing* if either of the following equivalent conditions holds (see [BR]):

(i) The representation $(\pi, V)$ does not contain finite dimensional subrepresentations.

(ii) The representation $(\pi \otimes \pi^*, V \otimes V^*)$ does not contain invariant vectors (here $(\pi^*, V^*)$ denotes the dual, or the contragradient, to $(\pi, V)$ representation).

(iii) For any unitary representation $(\sigma, W)$ the representation $(\pi \otimes \sigma, V \otimes W)$ does not contain fixed vectors.

LEMMA 8.2. *Let* $(\Omega, m)$ *be an ergodic* ME *coupling of* $\Gamma$ *with* $\Lambda$ *and let* $(\pi, V)$ *be some weakly mixing unitary* $\Lambda$-representation. *Then the the induced* $\Gamma$-representation $(\tilde{\pi}, \tilde{V})$ *has no fixed vectors.*

*Proof.* Assume that $\phi \in \tilde{V}$ is a $\tilde{\pi}(\Gamma)$-invariant unit vector; i.e. $\pi((\alpha(\gamma, x)^{-1}) \phi(\gamma \cdot x) = \phi(x)$ for all $\gamma \in \Gamma$. Define $F : \Omega \to V$ by $F(x\lambda) = \pi(\lambda^{-1}) \phi(x)$ for $x \in X$, $\lambda \in \Lambda$. Then for $x \in X$,

$$F(\gamma x) = F((\gamma \cdot x) \alpha(\gamma, x)) = \pi(\alpha(\gamma, x)^{-1}) \phi(\gamma \cdot x) = \phi(x) = F(x)$$

and $F(\gamma \omega) = F(\omega)$ for general $\omega \in \Omega$. Consider the restriction $f = F|_Y : Y \to V$. We have

$$f(y \cdot \lambda) = F(\beta(y, \lambda)^{-1} y \lambda) = F(y \lambda) = \pi(\lambda^{-1}) f(y)$$

and

$$\|f\|^2 = \int_Y \|F(y)\|^2 \, dm(y) = \frac{m(Y)}{m(X)} \int_X \|F(x)\|^2 \, dm(x) = \frac{m(Y)}{m(X)} \cdot \|\phi\|^2.$$

Now observe that the $\Lambda$-representation on the space $L^2(g : Y \to V)$ defined by

$$\lambda : g(y) \mapsto \pi(\lambda) g(y \cdot \lambda^{-1})$$

is isomorphic to the tensor product $(\pi \otimes \sigma, V \otimes L^2(Y))$ where $\sigma$ is the $\Lambda$-representation on $L^2(Y)$. The fact that $f \in V \otimes L^2(Y)$ is a fixed vector contradicts the assumption that $(\pi, V)$ is weakly mixing. $\square$

Recall that a countable (or more generally locally compact) group $\Gamma$ is said to have Kazhdan's property T if any unitary representation with almost invariant vectors necessarily contains an invariant vector. Corollary 1.4, asserting that property T is a ME invariant, can now be easily deduced from the following result of Bekka and Valette:

THEOREM 8.3 ([BV]). *If* $\Lambda$ *does not have Kazhdan's property* T, *then it admits a weakly mixing unitary representation* $(\pi, V)$ *with almost invariant vectors.*

Suppose that $\Gamma$ and $\Lambda$ are ME and $\Lambda$ does not have property T. Choose an ergodic ME coupling $(\Omega, m)$ of $\Gamma$ with $\Lambda$ (Lemma 2.2), and let $(\tilde\pi, \tilde V)$ be the $\Gamma$ representation induced from a weakly mixing $\Lambda$-representation $(\pi, V)$ which has almost invariant vectors. Then Lemmas 8.1 and 8.2 imply that the induced $\Gamma$-representation $(\tilde\pi, \tilde V)$ has almost invariant vectors, but does not have an invariant one.

*Remark* 8.4. A statement similar to Corollary 1.4 appears in [Zi3, Th.9.1.7 (b)], where the weak mixing assumption was imposed on the orbit equivalent group actions on the space. The result of Bekka and Valette enables us to bypass this assumption.

University of Illinois at Chicago, Chicago IL
*E-mail address*: furman@math.uic.edu